\documentclass[12pt]{amsart}

\usepackage{amssymb}

\usepackage{enumerate}

\usepackage{graphicx}

\makeatletter
\@namedef{subjclassname@2010}{%
  \textup{2010} Mathematics Subject Classification}
\makeatother

\theoremstyle{definition}

\numberwithin{equation}{section}

\begin{document}
\sloppy

\baselineskip=17pt



\title[On semi-barrelled spaces]{On semi-barrelled spaces}

\author[H. Bourl\`{e}s]{Henri Bourl\`{e}s}
\address{Satie, ENS de Cachan/CNAM \\61 Avenue Pr\'{e}sident Wilson\\ F-94230 Cachan,
France}
\email{henri.bourles@satie.ens-cachan.fr}

\date{}

\begin{abstract}
The aim of this paper is to clarify the properties of semi-barrelled spaces (also called countably
quasi-barrelled spaces in the literature).  These spaces were studied by several authors, in particular
in the classical book of N. Bourbaki "Espaces vectoriels
topologiques". However, six incorrect statements can be found in this reference. In particular: a Hausdorff and
quasi-complete semi-barrelled space is complete, a semi-barrelled,
semi-reflexive space is complete, a locally convex hull of semi-barrelled
semi-reflexive spaces is semi-reflexive, a locally convex hull of
semi-barrelled reflexive spaces is reflexive. \ We show through counterexamples that these
statements are false.  To conclude,  we show how these false claims can be corrected and
we collect some properties of semi-barrelled
spaces.
\end{abstract}



\maketitle

\section{Introduction}
Semi-barrelled spaces are studied in the classical book of N. Bourbaki \cite{Bourbaki} ("Espaces vectoriels
topologiques"), in  \S IV.3, n$%
{{}^\circ}%
$1 ("Espaces semi-tonnel\'{e}s"). \ A locally convex space $E$ is said to be semi-barrelled if the
following condition holds: (a) let $U$ be a bornivorous part of $E$ which is
the intersection of a sequence of convex balanced closed neighborhoods of $0$
in $E;$ then $U$ is a neighborhood of $0$ in $E.$ \ Chronologically, ($%
\mathcal{DF}$) spaces were introduced before semi-barrelled spaces, by
Grothendieck \cite{Grothendieck-F-DF}: a locally convex space $E$ is a ($%
\mathcal{DF}$) space if (a) holds and (b) the canonical bornology of $E$ has
a countable base. \ Then Husain \cite{Husain} considered spaces satisfying
(a) but not (b) and called them "countably semi-barrelled spaces". \ Thus,
these spaces are called "semi-barrelled spaces" by Bourbaki \cite%
{Bourbaki}, and we adopt the latter terminology in the sequel.

Bourbaki's account on semi-barrelled spaces consists of the above-quoted section, where three equivalent definitions are given, and of nine statements
(Exerc. 9, p. IV.60).  Six of them are not correct, therefore  a clarification of the properties of these spaces is needed and presented below.  These statements are the following ones:

(1) Let $E$ be a locally convex Hausdorff semi-barrelled space, $M$ be a
closed subspace of $E,$ $E^{\prime }$ be the dual of $E.$ \ The strong
topology $\beta \left( M^{0},E/M\right) $ is identical to the topology
induced on $M^{0}$ by the strong topology $\beta \left( E^{\prime },E\right)
.$

(2) Let $E$ be a locally convex Hausdorff space, $M$ a (non necessarily
closed) vector subspace of $E$. \ If $M$ is a semi-barrelled space, then the
strong topology $\beta \left( E^{\prime }/M^{0},M\right) $ is identical to
the quotient topology by $M^{0}$ of the strong topology $\beta \left(
E^{\prime },E\right) .$

(3) A Hausdorff and quasi-complete semi-barrelled space $M$ is complete.

(4) Let $E$ be a semi-barrelled Hausdorff space and $M$ be a closed subspace
of $E.$ \ Then $E/M$ is a semi-barrelled space.

(5) Let $\left( E_{n}\right) $ be a sequence of semi-barrelled spaces, $E$
be a vector space, and for each $n,$ let $f_{n}$ be a linear mapping from $%
E_{n}$ into $E.$ \ Suppose that $E$ is the union of the $f_{n}\left(
E_{n}\right) $. \ Then $E$ is semi-barrelled for the finest locally convex
topology for which all the $f_{n}$ are continuous.

(6) The completion of a semi-barrelled Hausdorff space is semi-barrelled.

(7) A semi-barrelled, semi-reflexive space $M$ is complete.

(8) Let $\left( E_{n}\right) $ and $E$ be as in (5). \ If each $E_{n}$ is
semi-reflexive and if $E$ is Hausdorff, then $E$ is semi-reflexive.

(9) Let $\left( E_{n}\right) $ and $E$ be as in (5). \ If each $E_{n}$ is
reflexive and if $E$ is Hausdorff, then $E$ is reflexive.

Claims (4) and (5) were proved in (\cite{Husain}, Thm. 8 and Corol. 14)
with the sequences $\left( E_{n}\right) $ and  $\left(f_{n}\right) $ replaced by 
non-necessarily countable families $\left( E_{\alpha}\right) $ and  
$\left(f_{\alpha}\right) $ in Claim (5) --and one can notice
 that it is not necessary to assume in Claim (4) that $M$ is
closed. \ Claim (6) was proved in (\cite{Husain-Khaleelulla}, Chap. V, Prop.
3, p. 133). \ We show below through counterexamples that Claims (1)-(3) and
(7)-(9) are false. \ So, Bourbaki's account on semi-barrelled spaces is very
misleading for people studying topological vector spaces, and a clarification of the
properties of these spaces is needed. \ We conclude by
additional remarks where we show how these false claims can be corrected and
we collect some properties of semi-barrelled
spaces.

\section{Counterexamples}
\subsection{Counterexample to Claim (1)}

An example is given in (\cite{Kothe}, \S 27, 2., p. 370) of a Montel space\ $%
F$ and a closed subspace $H$ of $F$, such that the initial topology of $F$
does not induce on $H$ its strong topology. \ So, let $E$ be the strong dual
of $F$ and $M=H^{0}.$ \ Then $E$ is a Montel space, thus it is barrelled, $M$
is a closed subspace of $E,$ $H=M^{0}$ according to the bipolar theorem, and
the strong topology $\beta \left( M^{0},E/M\right) $ is not identical to the
topology induced on $M^{0}$ by the strong topology $\beta \left( E^{\prime
},E\right) .$

\subsection{Counterexample to Claim (2)}

It is well-known that there exist non-distinguished
Fr\'{e}chet spaces (\cite{Kothe}, \S 31, 7., p. 435).\ Such a space $M$ can be embedded in the product $E$ of a
family of Banach spaces (\cite{Kothe}, \S 18, 3., (7), p. 208). \ Then the
strong dual $E_{\beta }^{\prime }$ is a locally convex direct sum of Banach
spaces, thus is barrelled; therefore $E_{\beta }^{\prime }/M^{0}$ (the
topology $\mathfrak{T}_{1}$ of which is the quotient topology of $\beta
\left( E^{\prime },E\right) $ by $M^{0}$)\ is barrelled. \ Since $M$ is
non-distinguished, its strong dual (the topology $\mathfrak{T}_{2}$ of which
is $\beta \left( E^{\prime }/M^{0},M\right) $) is not barrelled. \ Therefore
the topologies $\mathfrak{T}_{1}$ and $\mathfrak{T}_{2}$\ are not identical.

\subsection{Counterexample to Claims (3) and (7)}

An example of a non-complete Montel space has been given in (\cite{Komura}%
, \S 5). \ This space is barrelled and reflexive, hence quasi-complete (\cite%
{Bourbaki}, p. IV.16, Remarque 2). \ Therefore, there exist reflexive (thus
quasi-complete) barrelled spaces which are not complete.

\subsection{Counterexample to Claims (8) and (9)}

If Claim (8) is correct, so is Claim (9) since $E$ is barrelled if the
spaces $E_{n}$ are barrelled (\cite{Kothe}, \S 27, 1., (3), p. 368). \
However, in (\cite{Kothe}, \S 31, 5., p. 434), an example is given of a Fr%
\'{e}chet-Montel space $E_{1}=\lambda $ and a closed subspace $N$ of $E_{1}$
such that $E=E_{1}/N$ is topologically isomorphic to $l^{1},$ thus is not
reflexive. \ The space $E_{1}$ is reflexive and $E=f_{1}\left( E_{1}\right)
, $ where $f_{1}$ is the canonical surjection, is not. \ Therefore, Claim
(9) is false and Claim (8) is false too.

\section{Concluding remarks}

According to the proof of (\cite{Grothendieck-F-DF}, corol., p. 79), Claims
(8) and (9) are correct if $E$ and the sequence $\left( E_{n}\right) $ are
as in (5), assuming that every bounded part of $E$ is included in the closed
balanced convex hull of a finite number of $f_{n}\left( B_{n}\right) $ where
each $B_{n}$ is a bounded part of $E_{n}$ (in particular, this holds if
$E$ is a regular inductive limit of the $E_{n}$  \cite%
{Floret}, and  \textit{a fortiori}  if the $%
E_{n}$ are closed subspaces of $E$ and $E$ is the strict inductive limit of
the $E_{n}$). \ On the other hand, all statements are correct if ($\mathcal{%
DF}$) spaces are considered in place of semi-barrelled spaces (\cite%
{Grothendieck-F-DF}, prop. 5, p. 76 and corol., p. 79), (\cite%
{Grothendieck-book}, corol. 2, p. 170).   

A closed subspace $F$
of a semi-barrelled  space $E$ is not necessarily
semi-barrelled, as shown in (\cite%
{Grothendieck-F-DF}, p. 97) and (\cite{Iyahen}, (iii)), except if $F$ is finite-codimensional (\cite{Webb}, Thm. 6, p. 169) or if $F$ 
is countable-codimensional and such that for every bounded subset $B$ of $E$, $F$ is finite-codimensional in the space spanned by 
$F\cup B$  (\cite{Webb 73}, Thm. 3). \ Every infra-barrelled (also called \emph{quasi-barrelled}) space
is semi-barrelled (\cite{Husain}, Prop. 2, p. 292), every infra-barrelled
space is a Mackey space (\cite{Grothendieck-book}, p. 107), but a
semi-barrelled space is not necessarily a Mackey space by (\cite%
{Grothendieck-F-DF}, Remarque 8, p. 74) -- since neither is a ($\mathcal{DF}$%
) space. \ Further results on semi-barrelled spaces can be found in, e.g., 
\cite{Rademovic}, \cite{Krassovska}, \cite{Katsaras}, and (\cite{Khaleelulla}%
, Chap. III, p. 33 and n$%
{{}^\circ}%
$41., 44.; Chap. IV, 16.; Chap. VI, 4.).

\end{document}